\documentclass[a4paper, 10pt]{amsart}

\usepackage[numeric, lite, initials]{amsrefs}

\usepackage{amssymb}
\usepackage{amsmath}
\usepackage{mathrsfs}

\input xy
\xyoption{arrow} \xyoption{matrix}

\numberwithin{equation}{section}

\newcommand{\Z}{{\mathbb Z}}

\def\:{\colon}
\def\.{\cdot}
\def\lra{\longrightarrow}
\def\smash{\wedge}
\def\epsilon{\varepsilon}
\def\phi{\varphi}
\def\subset{\subseteq}

\def\geq{\geqslant}
\def\D{\mathscr D}

\def\o{\circ}
\def\dr{\mathscr D_R}
\def\de{\mathscr D_E}

\DeclareMathOperator {\Tor}{Tor} \DeclareMathOperator{\gr}{gr}
\DeclareMathOperator{\Sym}{Sym} 
\DeclareMathOperator{\Hom}{Hom}

\DeclareMathOperator{\coker}{coker}
\DeclareMathOperator{\holim}{holim}

\newtheorem{theorem}{Theorem}
\newtheorem{lemma}[theorem]{Lemma}
\newtheorem{proposition}[theorem]{Proposition}
\newtheorem{corollary}[theorem]{Corollary}

\numberwithin{equation}{section}

\theoremstyle{definition}
\newtheorem*{conventions}{Conventions}
\newtheorem{remark}[theorem]{Remark}

\title{Homology of $I$--adic towers}
\author{Samuel W\"{u}thrich}
\date{29/11/2005}


\address{Department of Pure Mathematics, University of Sheffield,
Hicks Building, Hounsfield Road, Sheffield S3 7RH, United
Kingdom.}

\email{S.Wuethrich@sheffield.ac.uk}

\begin{document}

\begin{abstract}
Let $R$ be a commutative ring and $I\lhd R$ an ideal generated by
a regular sequence. Then it is known that the natural sequences
\[
0 \to \Tor_*^R(R/I, I^s) \to \Tor_*^R(R/I, I^s/I^{s+1}) \to
\Tor_{*-1}^R(R/I, I^{s+1}) \to 0
\]
are short exact sequences of graded free $R/I$--modules, for any
$s\geq 0$. The aim of this paper is to give a proof which accounts
for the structural simplicity of the statement. It relies on a
minimum of technicalities and exposes the phenomenon in a
transparent way as a consequence of the regularity assumption. The
ideas discussed here are used in \cite{W2} to obtain a better
qualitative understanding of $I$-adic towers in algebraic
topology.
\end{abstract}

\maketitle

\section*{Introduction}

Let $R$ be a commutative ring with unit and assume that $I\lhd R$
is an ideal generated by a finite regular sequence. Consider the
$I$-adic filtration of $R$
\[
\cdots \subseteq I^{s+1} \subseteq I^s \subseteq \cdots \subseteq
I\subseteq R.
\]

\begin{theorem}\label{thm}
The short exact sequences
\begin{equation}\label{defses} 0 \lra
I^{s+1} \lra I^s \lra I^s/I^{s+1} \lra 0
\end{equation}
induce short exact sequences of graded free $R/I$-modules
\begin{equation}\label{ses}
0 \to \Tor_*^R(R/I, I^s) \to \Tor_*^R(R/I, I^s/I^{s+1}) \to
\Tor_{*-1}^R(R/I, I^{s+1}) \to 0.
\end{equation}
\end{theorem}

This result appears to have been known to commutative algebraists
for a while, but the author knows of no convenient self-contained
account in the algebraic literature. He encountered the statement
in Baker's preprint \cite{B}, which was written in view of
applications to algebraic topology. A published proof of Theorem
\ref{thm} based on Baker's ideas can be found in \cite{W}.

In this paper, we investigate the situation from a new point of
view. One might say that the moral content of Theorem \ref{thm} is
that the powers of the ideal $I$, which by birth are firmly rooted
in the world of $R$-modules, are in the given situation in fact
very neatly organized with respect to each other in the {\em
derived} category $\dr$ of $R$. This is made precise in Remark
\ref{rem}. We aim to accommodate this view by giving a proof which
exposes the statement in a conceptually simple way as a
consequence of the regularity assumption and which uses only a
minimum of technicalities. In particular, we completely avoid the
use of explicit chain complexes other than Koszul complexes. The
motivation for our work comes from algebraic topology, as we
briefly indicate now.

The basic objects of study in stable homotopy theory are ring
spectra. They represent cohomology theories which are equipped
with a multiplicative structure. One of the great achievements of
modern stable homotopy theory is the construction of a derived
category $\de$ of module spectra over $E$ for highly structured
ring spectra $E$ (see \cites{EKMM, HSS} for two different
accounts). This category shares many formal similarities with the
derived category of an ordinary ring. Therefore one hopes that an
understanding of algebraic results such as Theorem \ref{thm} in
derived terms might lead to an insight into derived categories of
ring spectra. That this is the case for the situation we study
here is shown in \cite{W2}.

To explain our approach to the proof of Theorem \ref{thm}, we
consider the diagram
\begin{equation*}\begin{array}{r}\label{dia}
\xymatrix@C=-0.8cm@R=0.6cm{0 \ar[rr] && S \ar[rr] && \Tor_*^R(S,
S) \ar[rr]^-{\delta^0}\ar[dr] && \Tor_{*-1}^R(S, I/I^2)
\ar[rr]^-{\delta^1}\ar[dr] &&
\Tor_{*-2}^R(S, I^2/I^3) \ar[rr]^-{\delta^2}\ar[dr] && \cdots\\
&\makebox[1.8cm]{} &&\makebox[1.8cm]{} &&\Tor_{*-1}^R(S, I)\ar[ru]
&& \Tor_{*-2}^R(S, I^2) \ar[ru] &&\makebox[1.8cm]{} && \cdots
 }
\end{array}
\end{equation*}
which is obtained by pasting together the sequences \eqref{ses}.
Here and in the remainder of this section, $S$ denotes $R/I$.
Clearly, the horizontal sequence in the diagram is exact if all
the sequences \eqref{ses} are exact. A first simple but important
observation is that the converse also holds, as a consequence of
the fact that short exact sequences give rise to long exact
sequences on $\Tor_*^R(S,-)$. Now the regularity assumption
implies that the $S$-modules $I^s/I^{s+1}$ are free and that
$\Tor_*^R(S,S)$ is isomorphic to a graded exterior algebra
$\Lambda_{S}(e_1, \ldots, e_n)$ on generators $e_i$ of degree
\nolinebreak $1$. We denote this exterior algebra by $\Lambda_*$
for short. Thus our sequence is of the form
\begin{equation}\label{modelles}
0 \lra S \xrightarrow{\  \epsilon\ } \Lambda_* \xrightarrow{\
\delta^0\ } I/I^2 \otimes_S \Lambda_{*-1}\xrightarrow{\  \delta^1\
} I^2/I^3 \otimes_S \Lambda_{*-2} \xrightarrow{\  \delta^2\ }
\cdots,
\end{equation}
where $\epsilon$ is the inclusion. By a naturality argument, we
can identify the maps $\delta^s$ as the connecting homomorphisms
associated to the short exact sequences
\[
0 \lra I^{s+1}/I^{s+2} \lra I^s/I^{s+2} \lra I^s/I^{s+1} \lra 0.
\]
Now the crucial observation is that $\delta^0$ determines all the
other maps in the sequence in a simple manner. This follows from
regarding the maps $\delta^s$ as the components of the connecting
homomorphism $\delta^*$ associated to the single short exact
sequence
\begin{equation}\label{singexp}
0 \lra \bigoplus_{s\geq 0} I^{s+1}/I^{s+2} \lra \bigoplus_{s\geq
0} I^s/I^{s+2} \lra \bigoplus_{s\geq 0} I^s/I^{s+1} \lra 0
\end{equation}
and noting that this is in fact a singular extension of algebras.
For the connecting homomorphism of such the following statement
holds, which we guess is well-known, but are unable to find in the
literature.
\begin{proposition}\label{prop}
Assume we are given a singular extension of $R$-algebras
\begin{equation*}\label{sing} 0 \lra J \lra A \lra
A' \lra 0
\end{equation*}
and some $R$-algebra $T$. Then the connecting homomorphism
\[
\partial\: \Tor_*^R(T, A') \lra \Tor_{*-1}^R(T, J)
\]
is a derivation with respect to the natural biaction of
$\Tor_*^R(T, A')$ on $\Tor_*^R(T,J)$.
\end{proposition}
Applied to \eqref{singexp}, this shows that $\delta^*$ is
determined by $\delta^1$, as follows. Regularity of the ideal $I$
implies that $\gr_I^*(R) = \bigoplus_{s\geq 0} I^s/I^{s+1}$ is
isomorphic to the symmetric algebra on $I/I^2$. Therefore
$\Tor_*^R(S,\gr_I^*(R))$ is the free graded commutative algebra
over $\Tor_*^R(S,S)$ generated by $\Tor_*^R(S, I/I^2)$. So we are
left to identify $\delta^1$ and to show that the sequence
\eqref{modelles} that it determines is exact. We do this by
identifying it with a well-known relative injective resolution of
$S$ over the coalgebra $\Lambda_*$, which we call the model
complex. It is obtained by rearranging the Koszul complex for the
regular sequence $(- x_1 + I^2, \ldots, - x_n + I^2)$ on
$\gr_I^*(R)$. We also explain at this point what the structure on
the Koszul complex is which makes \eqref{modelles} a sequence of
$\Lambda_*$-comodules. The statement in Theorem \ref{thm} about
the freeness of the $\Tor$-groups over $S$ is an easy consequence
of the identification of \eqref{modelles} as the model complex.

The paper is organized as follows. In Section \ref{sectionkoszul},
we recall basic material on Koszul complexes and regular
sequences. In Section \ref{sectionmodel}, we discuss the model
complex mentioned above. In Section \ref{sectionextension}, we
prove Proposition \ref{prop}. In Section \ref{sectioniadic}, we
identify the sequence \eqref{modelles} as the model complex
(Proposition \ref{sequence}) and deduce Theorem \ref{thm}. We
finish by explaining how the theorem leads to a characterization
of the $I$-adic filtration in the derived category $\dr$.

Many thanks go to Alain Jeanneret and Andrew Baker for their
support and guidance. Andrew Baker's paper \cite{B} was very
inspiring to me. I am grateful to the Swiss National Science
Foundation for financial support.

\begin{conventions}
All rings are assumed to be commutative and to have a unit. All
algebras are required to have a unit. We use the convention
$M^*=M_{-*}$ for graded objects and write $1$ for all identity
maps. An unlabelled $\otimes$ means $\otimes_R$, where $R$ is some
ring which is specified in the context.
\end{conventions}

\section{Koszul complexes and regular sequences}\label{sectionkoszul}

Let $(r_1, \ldots, r_n)$ be a sequence of elements of a ring $R$,
generating an ideal $I\subset R$. Consider the exterior algebra
$\Lambda_*^R(M)$ on an $R$--module $M$. It is naturally graded,
with $M$ as the submodule of degree one. We denote the degree of a
homogeneous element $e\in\Lambda_*^R(M)$ by $|e|$ and the product
on $\Lambda_*^R(M)$ by $\smash$. Recall that $\Lambda^R_*(M)$ is a
graded bialgebra, i.e.\@ it supports a coproduct
\begin{equation}\label{coproduct}
\Delta\: \Lambda^R_*(M)\lra \Lambda_*^R(M)\otimes\Lambda^R_*(M)
\end{equation}
which is a map of graded algebras. On a homogeneous element $m$ of
degree one, it is defined by $\Delta(m)=1\otimes m + m\otimes 1$.
The coproduct is graded cocommutative, in the sense that
$\Delta=\tau\o\Delta$. Here $\tau$ denotes the twist map, defined
by $\tau(e\otimes f) = (-1)^{|e| |f|} f\otimes e$ for two
homogeneous elements $e$ and $f$.

Now assume that $U$ is a free $R$--module of rank $n$ with fixed
basis $e_1, \ldots, e_n$. The $R$--linear map $d_1\: U\to R$ which
sends $e_j$ to $r_j$ has a unique extension to a graded
$R$--derivation
\[
d_*\: \Lambda^R_*(U) \lra \Lambda^R_{*-1}(U).
\]
So the image of a product of two homogeneous elements $e$ and $f$
is given by
\begin{equation}\label{derivationformula}
d_*(e\smash f) = d_*(e)\smash f + (-1)^{|e|} e \smash d_*(f).
\end{equation}
Formally, we may interpret $d_*$ as the expression
\[
d_* = \sum_{j=1}^n r_j \frac{\partial}{\partial e_j}.
\]
It follows easily that $d_*$ is a differential. To verify this, it
suffices to show that $(d_*\o d_*)(g)=0$ for all homogeneous $g$.
We prove this by induction over $|g|$. The case $|g|=1$ is
trivial. For $|g|>1$, write $g$ as a linear combination of
elements of the form $e\smash f$, with $e$ and $f$ homogeneous and
$|e|=1$. Then the inductive step follows from
\eqref{derivationformula}.

We do not distinguish between the concepts of differential graded
modules and chain complexes. In particular, we regard
$\Lambda_*^R(U)$ as a chain complex.

Another formulation of the fact that $d_*$ is a derivation is the
following. Let $(\Lambda^R_*(U)\otimes \Lambda^R_*(U),
d_*^\otimes)$ be the tensor product of the chain complex
$(\Lambda^R_*(U), d_*)$ with itself. Its components are given by
\[
(\Lambda^R_*(U)\otimes \Lambda^R_*(U))_n = \bigoplus_{p+q=n}
\Lambda^R_p(U)\otimes \Lambda^R_q(U).
\]
The value of its differential $d_*^\otimes$ on an element of the
form $e\otimes f$, where $e$ and $f$ are homogeneous, is given by
\[
d_*^\otimes(e\otimes f) = d_*(e)\otimes f + (-1)^{|e|} e\otimes
d_*(f).
\]
The reader may check that we can write
\begin{equation}\label{diffotimes}
d_*^\otimes = d_*\otimes 1 + \tau\o(d_*\otimes 1)\o\tau.
\end{equation}
To say that $d_*$ is a derivation is equivalent to state that the
product
\[
\smash\: \Lambda^R_*(U)\otimes\Lambda^R_*(U)\lra \Lambda^R_*(U)
\]
is a map of chain complexes. The coproduct $\Delta$, on the other
hand, is not compatible with the differential in such a way.
Instead, we have
\begin{equation}\label{coprodrep}
d_*^\otimes\o \Delta = 2(\Delta \o d_*).
\end{equation}
This is a consequence of the formula
\begin{equation}\label{coprodcol}
\Delta\o d_* = (d_*\otimes 1) \o \Delta,
\end{equation}
which we prove in an instant. Namely, cocommutativity of
$\Lambda_*^R(U)$ and \eqref{coprodcol} imply
\[
\tau\o(d_*\otimes 1)\o\tau\o \Delta= \tau\o (d_*\otimes 1)\o\Delta
= \tau\o\Delta\o d_* = \Delta\o d_*.
\]
Applying \eqref{coprodcol} again, \eqref{coprodrep} follows from
\eqref{diffotimes}.

So let us prove \eqref{coprodcol}. We abbreviate $\Lambda_*^R(U)$
by $\Lambda_*$ in the following. Regard $\Lambda_*\otimes
\Lambda_*$ as a $\Lambda_*$--bimodule by restricting scalars along
$\Delta$. So the left (right) action of an element $x\in\Lambda_*$
is given by left (right) multiplication by $\Delta(x)$. We claim
that both sides of \eqref{coprodcol} are derivations $\Lambda_*
\to (\Lambda_*\otimes\Lambda_*)_{*-1}$ with respect to this
bimodule structure. For the left hand side, note that as an
algebra map, $\Delta$ is a map of $\Lambda_*$--bimodules. As $d_*$
is a derivation, this implies by naturality that $\Delta\o d_*$ is
a derivation. For the right hand side, we use naturality in the
first argument. Namely, it can be checked that $d_*\otimes 1$ is a
derivation with respect to the canonical
$\Lambda_*\otimes\Lambda_*$--bimodule structure on
$\Lambda_*\otimes\Lambda_*$. Also, $\Delta$ is an algebra map.
Hence $(d_*\otimes 1)\o\Delta$ is a derivation with respect to the
$\Lambda_*$--bimodule structure obtained by restricting scalars
along $\Delta$. Now it is easily checked that both sides of
\eqref{coprodcol} coincide on $e_j$ and hence on all elements of
degree one. As these generate $\Lambda_*$, the statement follows.

The complex $(\Lambda^R_*(U), d_*)$ is called the Koszul complex
associated to the sequence $(r_1, \ldots, r_n)$; we denote it by
$K(r_1, \dots, r_n)$. The projection $\epsilon\: \Lambda_0^R(U) =
R\to R/I$ defines an augmentation of $K(r_1, \dots, r_n)$ over
$R/I$.

\begin{proposition}[\cite{M}*{Th.\@ 16.5}]\label{koszul}
For a regular sequence $(r_1,\ldots, r_n)$, the augmented
differential graded algebra $K(r_1, \dots, r_n)$ defines an
$R$--free resolution of $R/I$.
\end{proposition}

As applying $R/I\otimes -$ to $K(r_1, \ldots, r_n)$ kills all the
differentials, this implies

\begin{corollary}\label{tor}
There is an isomorphism of graded algebras
\[
\Tor_*^R(R/I, R/I) \cong R/I \otimes \Lambda_*^R(U).
\]
\end{corollary}

Let $\gr^*_I(R)$ be the graded algebra associated to the $I$--adic
filtration of $R$. Its components are the $R/I$--modules $
\gr^s_I(R) = I^s/I^{s+1}$, where by convention $I^0=R$. Let
$\{r_j\}$ denote the residue class of $r_j$ in $I/I^2$. For a ring
$S$ and an $S$--module $M$, let $\Sym^*_S(M)$ be the symmetric
algebra on $M$. If $M$ is free of rank $n$, $\Sym_S^*(M)$ is
isomorphic to a polynomial ring over $S$ in $n$ variables. Just as
the exterior algebra on $M$, $\Sym^*_S(M)$ admits a graded
bialgebra structure. On homogeneous elements $m$ of degree one,
the coproduct is given by $\Delta(m)=1\otimes m + m\otimes 1$.

\begin{proposition}[\cite{M}*{Th.\@ 16.2}]\label{gr}
For a regular sequence $(r_1,\ldots, r_n)$, $I/I^2$ is freely
generated as an $R/I$--module by $\{r_1\}, \ldots, \{r_n\}$, and
there is an isomorphism of graded algebras
\[
\gr^*_I(R) \cong \Sym_{R/I}^*(I/I^2).
\]
\end{proposition}

Assume now that the ground ring $R$ is graded and that $M$ is a
graded $R$--module. Then $\Lambda^R_*(M)$ and $\Sym^*_R(M)$ are
bigraded $R$--bialgebras. If $r_1, \ldots, r_n$ are homogeneous
elements, the ideal $I=(r_1, \ldots, r_n)$ is graded and
$\gr_I^*(R)$ is a bigraded $R/I$--algebra.

\section{A certain exact sequence}\label{sectionmodel}

Let $V$ and $W$ be free modules of rank $n$ over a given ring $S$,
with bases $e_1, \ldots, e_n$ and $f_1, \ldots, f_n$ respectively.
The sequence $(f_1, \ldots, f_n)$ is regular on the graded
$S$--algebra $R=\Sym^*_S(W)$. It generates an ideal $I$ such that
$R/I\cong S$. The Koszul complex
\begin{equation}\label{modelcomplex}
K(f_1, \ldots, f_n) = \Lambda_*^R(R\otimes_S V)
\end{equation}
provides by Proposition \ref{koszul} an $R$--free resolution of
$S$. By the remark at the end of Section \ref{sectionkoszul},
$K(f_1, \ldots, f_n)$ is bigraded. Clearly, $K(f_1, \ldots, f_n)$
is obtained as a bigraded $S$--bialgebra from
$\Lambda_*=\Lambda_*^S(V)$ by extending scalars, i.e.
\begin{equation}\label{modelcomplex2}
K(f_1, \ldots, f_n) \cong R\otimes_S \Lambda_*.
\end{equation}
In particular, the coproduct $\Delta$ on $K(f_1, \ldots, f_n)$
corresponds to $1\otimes\Delta'$, where $\Delta'$ is the coproduct
on $\Lambda_*$. We denote the differential of $R\otimes_S
\Lambda_*$ corresponding under \eqref{modelcomplex2} to the
differential $d_*$ of $K(f_1, \ldots, f_n)$ by $d_*$ as well.

Schematically, the complex \eqref{modelcomplex} can be depicted as
\begin{equation}\label{bigdiagram}\begin{array}{c}
\xymatrix@C=0.3cm@R=0.3cm{ \cdots\\
R^0 \otimes_S \Lambda_3\ar[rd] & \cdots \\
R^0 \otimes_S \Lambda_2 \ar[rd] & R^1 \otimes_S \Lambda_2\ar[rd]&
\cdots
\\
R^0 \otimes_S \Lambda_1 \ar[rd] & R^1 \otimes_S \Lambda_1 \ar[rd]
& R^2 \otimes_S \Lambda_1\ar[rd]& \cdots
\\
R^0 \otimes_S \Lambda_0 & R^1 \otimes_S \Lambda_0 & R^2 \otimes_S
\Lambda_0 & R^3 \otimes_S \Lambda_0 & \cdots. }
\end{array}
\end{equation}
Note that $R^0=\Lambda_0=S$. Clearly, the complex
\eqref{modelcomplex} is a resolution of $S$ precisely because each
diagonal sequence of length at least two is exact. Grouping the
diagram by columns instead of rows and prefixing the inclusion
$\epsilon\: S\to \Lambda_*$, we obtain an exact sequence of
$S$--modules
\begin{equation}\label{modelsequence}
0 \lra S \xrightarrow{\   \epsilon\  } \Lambda_* \xrightarrow{\
\partial^0\  } R^1 \otimes_S  \Lambda_{*-1}\xrightarrow{\
\partial^1\  }  R^2 \otimes_S  \Lambda_{*-2}
\xrightarrow{\  \partial^2\ } \cdots.
\end{equation}
We call \eqref{modelsequence} the {\em model complex under $S$ of
rank $n$} and write its differential as
\[
\sum_{j=1}^n f_i \otimes \frac{\partial}{\partial e_i}.
\]
Here $f_i$ means multiplication by $f_i$.

We have more structure on \eqref{modelsequence}. Recall that an
extended right $\Lambda_*$--comodule is one of the form
$M\otimes_S \Lambda_*$, where $M$ is an $S$--module and where the
coaction is given by $1\otimes \Delta$. More generally, a relative
injective $\Lambda_*$--comodule is a retract of an extended one.
Now endow $S$ with the trivial coaction and the other terms of
\eqref{modelsequence} with the extended coaction. Then we claim
that the coaugmentation $\epsilon$ and the differentials
$\partial^j$ are $\Lambda_*$--colinear. This is clear for
$\epsilon$. For the $\partial^j$, it is a reformulation of
equation \eqref{coprodcol} for the coproduct $\Delta$ on $K_* =
K(f_1, \ldots, f_n)$. To see this, consider
\[
\xymatrix@C=0cm@R=.7cm{ K_* \ar[rr]^-{d_*} \ar[dd]^-\Delta
\ar[dr]^-{\cong} && K_{*-1} \ar'[d]^-\Delta[dd]\ar[dr]^-{\cong}
\\
& R\otimes_S \Lambda_* \ar[rr]^<<<<<<<<<<{d_*}
\ar[dd]^<<<<<<{1\otimes\Delta'} && R\otimes_S \Lambda_{*-1}
\ar[dd]^-{1\otimes\Delta'}
\\
K_*\otimes_{R} K_* \ar'[r]^-{d_*\otimes 1}[rr] \ar[dr]^-{\cong} &&
K_{*-1}\otimes_{R} K_* \ar[dr]^-{\cong}
\\
& R\otimes_S \Lambda_*\otimes_S\Lambda_* \ar[rr]^-{d_*\otimes 1}
&& R\otimes_S \Lambda_{*-1}\otimes_S\Lambda_*.}
\]
By the remark after \eqref{modelcomplex2}, the top, bottom and
lateral faces of the cube commute. The back face commutes by
equation \eqref{coprodcol}. This forces the front face to commute,
which implies that the $\partial^j$ are maps of
$\Lambda_*$--comodules. The following proposition summarizes our
observations.

\begin{proposition}
The model complex \eqref{modelsequence} is a complex of right
$\Lambda_*$--comodules. It provides a relative injective
resolution of the trivial $\Lambda_*$--comodule $S$.
\end{proposition}

\section{Tor of a singular extension of
algebras}\label{sectionextension}

Let $T$ be an algebra over a ring $R$. The functor $\Tor_*^R(T,-)$
maps $R$--modules to graded $T$--modules and $R$--algebras to
graded $T$--algebras. If $R$ is a graded ring and $T$ a graded
$R$--algebra, $\Tor_*^R(T, -)$ takes graded $R$--modules to
bigraded $T$--modules and graded $R$--algebras to bigraded
$T$--algebras.

Let $p\: A\to A'$ be a surjection of $R$--algebras. Denote the
kernel of $p$ by $J$ and the inclusion $J\to A$ by $i$. Assume
that the extension
\begin{equation}\label{extension}
0\lra J \xrightarrow{\ \ i\ \ } A \xrightarrow{\ \ p\ \ } A' \lra
0
\end{equation}
is singular, i.e.\@ that the multiplication on $A$ restricted to
$J$ is trivial. An equivalent condition is that the $A$--bimodule
structure on $J$ lifts to $A'$. This action induces a bimodule
structure of $\Tor_*^R(T, A')$ on $\Tor_*^R(T, J)$; the action
maps are given by
\begin{align*}
& \Tor_*^R(T, A') \otimes_T \Tor_*^R(T, J) \lra \Tor_*^R(T,
A'\otimes J) \lra \Tor_*^R(T, J)
\\
& \Tor_*^R(T, J) \otimes_T \Tor_*^R(T, A') \lra \Tor_*^R(T,
J\otimes A') \lra \Tor_*^R(T, J)
\end{align*}
where the first maps are K\"unneth maps and the second maps are
induced by the left and the right $A'$--actions on $J$.

We need some preparation for the proof of Proposition \ref{prop}.
First, recall the ``Fundamental Lemma'' from homological algebra.
Let $M$ and $N$ be $R$--modules, let $P_*\to M$ be a complex over
$M$ with projective components and let $Y_*\to N$ be a resolution
of $N$. Then a given map $f\: M\to N$ can be lifted to a chain map
$F\: P_*\to Y_*$ covering $f$, and this lift is unique up to
homotopy \cite{CE}*{Prop.\@ V.1.1}. Now according to the Horseshoe
Lemma \cite{CE}*{Prop.\@ V.2.2}, it is possible to choose
projective resolutions $P_*\to J$, $Q_*\to A$ and $Q'_*\to A'$ of
$J$, $A$ and $A'$ respectively which fit into a short exact
sequence of chain complexes
\[
\xymatrix{0\ar[r] & P_* \ar[r]^\iota\ar[d] & Q_* \ar[r]^\pi\ar[d]
& Q'_* \ar[d] \ar[r] & 0
\\ 0 \ar[r] & J\ar[r]^i & A \ar[r]^p & A' \ar[r] & 0.
}
\]
A tensor product $P\otimes Q$ of two projective modules $P$ and
$Q$ is projective, because $\Hom_R(P\otimes Q, -) = \Hom_R(P,
\Hom_R(Q, -))$ is exact. So the Fundamental Lemma implies that we
can lift the products on $A$ and $A'$ to chain maps
\[
\mu\: Q_* \otimes Q_* \lra Q_*, \quad \mu'\: Q'_*\otimes Q'_* \lra
Q'_*.
\]
This allows us to interpret $Q_*$ and $Q'_*$ as (not necessarily
associative) differential graded algebras. Similarly, we can
construct maps
\begin{align*}
& \gamma_l\: Q_*\otimes P_* \lra P_*,  & \gamma'_l\: Q'_*\otimes
P_* \lra P_*,
\\
& \gamma_r\: P_*\otimes Q_* \lra P_*, & \gamma'_r\: P_*\otimes
Q'_* \lra P_*,
\end{align*}
covering the left and right actions of $A$ and $A'$ on $J$,
respectively. As a consequence of the Fundamental Lemma, the
diagrams
\begin{equation}\label{firstdiagram}\begin{array}{c}
\xymatrix{ Q_*\otimes P_* \ar[r]^{1\otimes \iota}\ar[d]^{\gamma_l}
& Q_*\otimes Q_* \ar[r]^{\pi\otimes\pi} \ar[d]^{\mu}
& Q'_*\otimes Q'_* \ar[d]^{\mu'} \\
P_* \ar[r]^\iota & Q_* \ar[r]^\pi & Q'_*}
\end{array}
\end{equation}
and
\begin{equation}\label{seconddiagram}\begin{array}{c}
\xymatrix{ Q_*\otimes P_* \ar[r]^{\pi\otimes 1} \ar[d]^{\gamma_l}
& Q'_*\otimes P_* \ar[ld]^{\gamma'_l}
\\
P_*}
\end{array}
\end{equation}
commute in the homotopy category of chain complexes; similarly for
the analogous diagrams involving the right actions.

\begin{proof}[Proof of Proposition \ref{prop}]
Let $\alpha=\{x'\}\in \Tor_p^R(T, A')$ be the residue class of
some $x'\in T\otimes Q'_p$. Let $x$ be a lift of $x'$ to $T\otimes
Q_p$ and $v$ a lift of $d(x)$ to $T\otimes P_{p-1}$, where $d$ is
the differential of $T\otimes Q_*$. For $\beta \in \Tor_q^R(T,
A')$, choose elements $y'$, $y$ and $w$ in a similar way. By
definition, $\partial(\alpha)$ and $\partial(\beta)$ are
represented by $v$ and $w$ respectively. Because the diagram
\eqref{firstdiagram} above commutes, $\partial(\alpha \. \beta)$
is represented by a lift of $d(\mu(x\otimes y))$ to $T\otimes
P_{p+q-1}$. As $Q_*$ is a differential graded algebra, we have
\[
d(\mu(x\otimes y)) = \mu(d(x)\otimes y) + (-1)^p \mu(x\otimes
d(y)).
\]
From commutativity of \eqref{firstdiagram} and
\eqref{seconddiagram}, it follows that
\[
\{ \mu(x\otimes d(y))\} = \{\gamma_l(x\otimes w)\}=
\{\gamma'_l(x'\otimes w)\}= \alpha \. \partial(\beta)
\]
and similarly that $\{\mu(d(x)\otimes y)\} =
\partial(\alpha) \. \beta$. Altogether we have proved
\[
\partial(\alpha \. \beta) = \partial(\alpha)\. \beta
+ (-1)^p \alpha\. \partial(\beta),
\]
which was the claim.
\end{proof}

\section{The $I$--adic tower}\label{sectioniadic}

Let $(r_1, \ldots, r_n)$ be a regular sequence of a ring $R$,
generating an ideal $I$. Put $S=R/I$. Consider the short exact
sequences of $R$--modules
\begin{equation}\label{sesconnect}
\mathcal F^s\: \quad 0 \lra I^{s+1}/I^{s+2} \lra I^s/I^{s+2} \lra
I^s/I^{s+1} \lra 0
\end{equation}
for $s\geq 0$, where the maps are the canonical injections and
projections. Associated to the $\mathcal F^s$ are connecting
homomorphisms
\begin{equation}\label{connpartial}
\delta^s\: \Tor_*^R(S, I^s/I^{s+1}) \lra \Tor_{*-1}^R(S,
I^{s+1}/I^{s+2}).
\end{equation}
We may view their direct sum as an endomorphism $\delta^*$ of the
bigraded $S$--module
\[
\Tor_*^R(S, \gr_I^*(R)) \cong \bigoplus_{s\geq 0} \Tor_*^R(S,
I^s/I^{s+1})
\]
of bidegree $(-1,1)$. Let $p$\/ be the projection $R\to S$ and
$\eta$ the composition
\[
\eta\: S \cong \Tor_*^R(S, R) \xrightarrow{\  p_*\  } \Tor_*^R(S,
S).
\]
Let $W$ be a free $S$--module of rank $n$ with fixed basis $f_1,
\dots, f_n$. Proposition \ref{gr} implies that we obtain an
isomorphism $I/I^2\cong W$ of $S$--modules by mapping the residue
class $\{r_j\}$ to $-f_j$. This induces an isomorphism
$\Sym^*_S(I/I^2)\cong\Sym^*_S(W)$ of algebras. Precomposing it
with the algebra isomorphism $\gr^*_I(R)\cong \Sym^*_S(I/I^2)$
from Proposition \ref{gr}, we obtain an algebra isomorphism
\[
\phi\: \gr_I^*(R) \cong \Sym^*_{S}(W).
\]
The reason for the minus sign will become clear soon. Combining
$\phi$ with the isomorphism from Corollary \ref{tor}, we obtain an
isomorphism of bigraded algebras
\begin{equation}\label{identify}
\psi_*\: \Tor_*^R(S, \gr^*_I(R)) \xrightarrow{\ \ \cong\ \ }
\Lambda_*^{\Sym_{S}^*(W)}(\Sym_{S}^*(W) \otimes U),
\end{equation}
where $U$ is a free $R$--module of rank $n$. We fix a basis $e_1,
\dots, e_n$ of $U$ and write $e_j$ for the image of $e_j$ in
$V=S\otimes U$ as well. Recall that we have defined a differential
on the right side of \eqref{identify} in \eqref{modelcomplex}.

\begin{proposition}\label{sequence}
The endomorphism $\delta^*$ of $\Tor_*^R(S, \gr_I^*(R))$
corresponds under the isomorphism $\psi_*$ from \eqref{identify}
to the differential $d_*$ of the Koszul complex. In other words,
the sequence of graded $S$--modules
\begin{equation}\label{lestor}
0 \lra S \xrightarrow{\  \eta\  } \Tor_*^R(S, S) \xrightarrow{\
\delta^0\ } \Tor_{*-1}^R(S, I/I^2) \xrightarrow{\ \delta^1\ }
\cdots
\end{equation}
is mapped under $\psi_*$ to the model complex
\eqref{modelsequence} under $S$. In particular, \eqref{lestor} is
exact.
\end{proposition}

\begin{proof}
We need to identify the maps $\partial^s$ from
\eqref{modelsequence} with the connecting homomorphism $\delta^s$
from \eqref{connpartial} under the isomorphism $\psi_*$. Recall
that the direct sum over the $\partial^s$ is the differential
$d_*$ of the Koszul complex
\[
\Lambda_*^{\Sym_{S}^*(W)}(\Sym_{S}^*(W)\otimes V).
\]
We have constructed $d_*$ as the unique graded derivation which
maps the elements $e_j$ to $f_j$. So it suffices to show that the
direct sum over the $\delta^s$ has the corresponding properties.
To see this, consider the singular extension of algebras
\[
0 \lra \bigoplus_{s\geq 0} I^{s+1}/I^{s+2} \lra \bigoplus_{s\geq
0} I^s/I^{s+2} \lra \bigoplus_{s\geq 0} I^s/I^{s+1} \lra 0
\]
over $\bigoplus_{s\geq 0} I^s/I^{s+1} = \gr^*_I(R)$. Its
connecting homomorphism
\[
\delta^*\: \Tor_*^R(S, \bigoplus_{s\geq 0} I^s/I^{s+1}) \lra
\Tor_{*-1}^R(S, \bigoplus_{s\geq 0} I^{s+1}/I^{s+2})
\]
is the direct sum of the $\delta^s$. By Proposition \ref{prop},
$\delta^*$ is a derivation. So we are left to show that $\delta^0$
maps the element $e_j\in\Tor_1^R(S, S)$ to $-\{r_j\}\in I/I^2
\cong\Tor_0^{S}(S, I/I^2)$, which corresponds to $f_j$ under
$\psi_*$. This is the content of the lemma below.
\end{proof}

The following lemma shows why we had to define $\psi_*$ in the way
we did.

\begin{lemma}\label{connlemma}
The connecting homomorphism
\[
\delta^0\: \Tor_*^R(S, S) \lra \Tor_{*-1}^R(S, I/I^2)
\]
maps the element $e_j$ to the residue class $-\{r_j\}$ of $-r_i$
in $I/I^2$.
\end{lemma}

\begin{proof}
We have free $R$--resolutions $\Lambda_*^R(U)\to S$ and
$I/I^2\otimes\Lambda_*^R(U)\to I/I^2$. By the Horseshoe Lemma, we
may construct a differential $d'_*$ and an augmentation
$\epsilon'$ on $(I/I^2 \oplus R)\otimes \Lambda_*^R(U)$ over
$R/I^2$ such that the diagram
\[
\xymatrix@R=0.7cm@C=0.5cm{& 0\ar[d] & 0\ar[d]&0\ar[d]\\
0\ar[r] & I/I^2\otimes\Lambda_*^R(U) \ar[r] \ar[d] & (I/I^2\oplus
R)\otimes\Lambda_*^R(U) \ar[r]\ar[d]^-{\epsilon'} &
\Lambda_*^R(U) \ar[d]\ar[r] & 0 \\
0\ar[r] & \ \ I/I^2\ \ \ar[r]\ar[d] & \ \  R/I^2\ \ \ar[r]\ar[d] &
\ \ S\ \ \ar[d] \ar[r]
&0\\
&0&0&0 }
\]
is a short exact sequence of acyclic complexes. We can set
$\epsilon'(\{r_j\})= \{r_j\}$, where $\{r_j\}$ denotes the residue
classes of $r_j$ both in $I/I^2$ and in $R/I^2$. Then we can
define
\[
d_1'\bigl((0,1)\otimes e_j\bigr) = (r_j, -\{r_j\}) \in R\oplus
I/I^2.
\]
This shows that $\delta^0(e_j)= -\{r_j\}$.
\end{proof}

The projections $p_{s+1}\: I^{s+1} \to I^{s+1}/I^{s+2}$ and
$I^s\to I^s/I^{s+2}$ induce a morphism
\begin{equation*}\label{sesfinal}
\xymatrix{\mathcal E^s\: \quad  0 \ar[r] & I^{s+1}
\ar[r]\ar[d]^-{p_{s+1}} & I^s \ar[r]^-{p_s} \ar[d] & I^s/I^{s+1}
\ar[r]\ar@{=}[d] &  0\\
\mathcal F^s\: \quad 0 \ar[r] &  I^{s+1}/I^{s+2} \ar[r] &
I^s/I^{s+2} \ar[r] & I^s/I^{s+1} \ar[r] &  0}
\end{equation*}
of short exact sequences. It exhibits $\mathcal F^s$ as the
pushout $(p_{s+1})_*(\mathcal E^s)$ of $\mathcal E^s$. So by
naturality, the connecting homomorphism $\delta^s$ of $\mathcal
F^s=(p_{s+1})_*(\mathcal E^s)$ factors into
\[
\Tor^R_*(S, I^s/I^{s+1}) \xrightarrow{\  \epsilon^s\  }
\Tor_{*-1}^R(S, I^{s+1}) \xrightarrow{(p_{s+1})_*} \Tor_{*-1}^R(S,
I^{s+1}/I^{s+2}),
\]
where $\epsilon^s$ is the connecting homomorphism associated to
$\mathcal E^s$.

\begin{proof}[Proof of Theorem \ref{thm}]
Consider the diagram of graded $S$--modules
\begin{equation*}\label{diagram}\begin{array}{c}
\xymatrix@C=-0.7cm@R=0.8cm{0 \ar[rr] &  & S \ar[rr]^-\eta & &
\Tor_*^R(S, S) \ar[rr]^-{\delta^0}\ar[rd]_{\epsilon^0} &&
    \Tor_{*-1}^R(S, I/I^2) \ar[rr]^-{\delta^1}\ar[rd]_{\epsilon^1} &&
\Tor_{*-2}^R(S, I^2/I^3) \ar[rr]^-{\delta^2} && \cdots.
\\
&\makebox[1.6cm]{} &&\makebox[1.6cm]{} && \Tor_{*-1}^R(S, I)
\ar[ru]_-{(p_1)_*} && \Tor_{*-2}^R(S, I^2) \ar[ru]_-{(p_2)_*}
&&\makebox[1.6cm]{} }
\end{array}
\end{equation*}
We know from Proposition \ref{sequence} that the row is exact. So
the claim is that the diagram is isomorphic to
\begin{equation*}\label{diagram2}
\xymatrix@C=-0.9cm@R=0.8cm{0 \ar[rr] &  & S \ar[rr]^-\eta & &
\Tor_*^R(S, S) \ar[rr]^-{\delta^0}\ar[rd] && \Tor_{*-1}^R(S,
I/I^2) \ar[rr]^-{\delta^1}\ar[rd] && \Tor_{*-2}^R(S, I^2/I^3)
\ar[rr]^-{\delta^2} && \cdots
\\
&\makebox[1.9cm]{} &&\makebox[1.9cm]{} &&
\coker\eta\cong\ker\delta^1 \ar[ru] &&
\coker\delta^0\cong\ker\delta^2 \ar[ru] &&\makebox[1.9cm]{} }
\end{equation*}
We prove the statement by induction as follows. We know that
$\eta$ is injective, hence $\Tor_{*-1}^R(S, I)\cong\coker\eta$. By
exactness, we have $\coker\eta\cong\ker\delta^1$. Hence the
morphism $(p_1)_*$ is injective. It follows that $\Tor_{*-2}^R(S,
I^2)$ is isomorphic to $\coker(p_1)_*=\coker(\delta^0)$, and so
on.

For the second claim, we have to show that the submodule
$\ker\delta^s\cong\Tor_*^R(S, I^s)$ of the free $S$--module
$\Tor_*^R(S, I^s/I^{s+1})$ is free. By Proposition \ref{sequence},
this is equivalent to the kernel of $\partial^s$ from the model
complex \eqref{modelsequence} over $S$ being $S$--free. If
$S\cong\Z$, this is automatic. For general $S$, note that the
model complex over $S$ can be obtained by applying
$S\otimes_{\Z}-$ to the model complex over $\Z$. Therefore,
\[
\ker\bigl(\partial^s\: \Sym^s_S(W)\otimes_S \Lambda_*^S(U) \to
\Sym^{s+1}_S(W)\otimes_S \Lambda_{*-1}^S(U)\bigr)
\]
is isomorphic to
\[
S\otimes_\Z \ker\bigl(\partial^s\: \Sym^s_\Z(W')\otimes_\Z
\Lambda_*^\Z(V') \to \Sym^{s+1}_\Z(W')\otimes_\Z
\Lambda_{*-1}^\Z(V')\bigr),
\]
where $V'$ and $W'$ are free $\Z$--modules of rank $n$. Therefore
$\ker(\partial^s)$ is $S$--free.
\end{proof}

\begin{remark}\label{rem}
We can derive from the theorem a characterization of the tower
\begin{equation}\label{adamstower}
\begin{array}{c}
\xymatrix@!C=.8cm{ \ldots \ar[r]& I^3 \ar[r]\ar[d] & I^2
\ar[r]\ar[d] & I \ar[r]\ar[d] & R\ar[d] \\
& I^3/I^4\ar[lu]|-\circ & I^2/I^3\ar[lu]|-\circ &
I/I^2\ar[lu]|-\circ & S \ar[lu]|-\circ
 }
\end{array}
\end{equation}
in the derived category $\D_R$ of the ring $R$ (an arrow with a
circle indicates a map of degree $-1$). We do not go into detail
here as the situation is entirely analogous to the one in topology
considered in \cite{W2}. What we have proved implies, in the
language of injective classes (see for instance \cite{C} or
\cite{W2}), that the tower is an Adams resolution of $R$ with
respect to the injective class associated to $S=R/I$. As such, it
is uniquely characterized up to isomorphism by the sequence
\begin{equation*}
0 \lra R \lra S \lra \Sigma I/I^2 \lra \Sigma^2 I^2/I^3 \lra
\cdots,
\end{equation*}
derived from \eqref{adamstower} in an obvious way. As
$I^s/I^{s+1}$ is a direct sum of suspended copies of $S$, we need
only $S$ and its endomorphisms in $\D_R$ to describe the sequence.

We also obtain a characterization of the completion $R\hat{_I}$ of
$R$ with respect to $I$ in $\mathscr{D}_R$. In our context, it
appears as the homotopy limit $\holim_s R/I^s$. We find that
$\holim_s R/I^s$ is the completion $R\hat{_S}$ of $R$ with respect
to $S$ in $\mathscr D_R$, in the sense of Dwyer and Greenlees
\cite{GD}. Namely, Prop.\@ 6.14 in this paper states that
completion with respect to $S$ is the same as Bousfield
localization for the homology theory $S_*(-)=H_*(S\otimes^{\bf L}
-)$, where $\otimes^{\bf L}$ is the derived tensor product. Now
all the $R/I^s$ are $S_*$--local, and Theorem \ref{thm} implies
that $R\to\holim_s R/I^s$ is an $S_*$--equivalence.
\end{remark}

\end{document}